\documentclass{ifacconf}

\usepackage{graphicx}      %
\usepackage{natbib}        %
\usepackage{amsmath}

\usepackage{tikz}
\usepackage{pgfplots}
\pgfplotsset{compat=1.11}
\usepackage{amssymb}
\usepackage{mathtools}

\usepackage[
  algo2e,
  commentsnumbered,
  linesnumbered,
  ruled
]{algorithm2e}
\SetKwProg{Fn}{Function}{\string:}{end}
\SetKwProg{Find}{find}{\string:}{}
\SetKw{st}{\ such that}

\newcommand{\TheTitle}{Exponential Convergence of Online Enrichment in Localized Reduced Basis Methods}

\newcommand{\norm}[1]{{\left\lVert{#1}\right\rVert}}
\newcommand{\spanset}{\operatorname{span}}
\newcommand{\argmax}{\operatornamewithlimits{arg\,max}}
\newcommand{\argmin}{\operatornamewithlimits{arg\,min}}
\newcommand{\localenrichment}{\widehat u_n}
\newcommand{\localenrichedsolution}[1]{\widehat u_{n,#1}}
\begin{document}
\begin{frontmatter}

\title{\TheTitle}

\author[First]{Andreas Buhr}

\address[First]{
  Institute for Computational and Applied Mathematics,
    University of M\"unster, Einsteinstra\ss e 62, 48149 M\"unster, Germany.
    (e-mail: andreas@andreasbuhr.de
).
}

\begin{abstract}                %
Online enrichment is the extension of a reduced solution space
based on the solution of the reduced model.
Procedures for online enrichment were published for many localized model order
reduction techniques.
We show that residual based online enrichment on overlapping domains
converges exponentially.
Furthermore, we present an optimal enrichment strategy which couples the
global reduced space with a local fine space.
Numerical experiments on the two dimensional stationary heat equation with high contrast and channels
confirm and illustrate the results.
\end{abstract}


\end{frontmatter}

\section{Introduction}
Online enrichment is the extension of an existing reduced solution space
based on the solution of the reduced model.
It is used in many localized model order reduction techniques:
ArbiLoMod, introduced in \cite{Buhr2017}, employs online enrichment.
For the Localized Reduced Basis Multiscale Method (LRBMS) which was introduced in \cite{Albrecht2012},
online enrichment was presented in \cite{Albrecht2013} and in \cite{Ohlberger2015}.
For the Generalized Multiscale Finite Element Method (GMsFEM) which was introduced in \cite{Efendiev2013},
online enrichment was discussed in \cite{Chung2015a}.
For the Constraint Energy Minimizing Generalized Multiscale Finite Element
Method (CEM-GMsFEM) which was introduced in \cite{Chung2017a}, 
online enrichment was discussed in \cite{Chung2017c}.

However, no a priori convergence analysis is available for the ArbiLoMod and the LRBMS.
For the GMsFEM and CEM-GMsFEM, exponential convergence was shown in \cite{Chung2015a} and \cite{Chung2017c},
but only for the case where the reduced space already has good properties.
The proofs given here do not require any properties of the reduced space.

The online enrichment presented here is usually
applied within a localized model order reduction method
like the methods mentioned above.
To analyze the convergence behavior and obtain
a priori estimates, we isolated the online enrichment.

While the setting presented in the following resembles
overlapping domain decomposition methods,
it is different as it generates reduced solution
spaces in each iteration. The approximate solution
in each step is the solution of the reduced problem,
in contrast to domain decomposition methods, were
the approximate solution is calculated as the sum
of the previous approximation and a correction term.
\section{Setting}
On the domain $\Omega$ with $\dim(\Omega) =d \in \{2,3\}$
we approximate the solution of the stationary heat equation
\begin{equation}
-\nabla(\kappa \nabla u) = f
\end{equation}
with homogeneous Dirichlet boundary conditions.
Non-homogeneous Dirichlet boundary conditions
could be handled with the usual shifting technique.
The methods presented below could easily
be extended for Neumann boundary conditions.
$\kappa$ is the heat conductivity.
In the space $V := H^1_0(\Omega)$, $u \in V$ 
is the unique solution of
\begin{equation}
a(u,\varphi) = f(\varphi) \quad \forall \varphi \in V
\end{equation}
with the coercive bilinear form 
\begin{equation}
a(u,\varphi) := \int_\Omega \kappa \nabla u \nabla \varphi dx
\end{equation}
and the linear form
\begin{equation}f(\varphi) := \int_\Omega f \varphi dx.
\end{equation}
We assume an overlapping domain decomposition with 
$N_D$ subdomains $\omega_i$ so that $\bigcup_{i=1,\dots, N_D} \omega_i = \Omega$.
We define local spaces $O_i := H^1_0(\omega_i)$
and assume a partition of unity
$\varrho_i \in O_i$,
$1 \equiv \sum_{i=1,\dots,N_D}\varrho_i$.
The partition of unity is not used in the algorithm,
it is only required for the proofs.
On the space $V$ and the spaces $O_i$, we use
the energy inner product and norm induced by $a$.

\section{Enrichment Algorithm}
\IncMargin{1em}
\begin{algorithm2e}
\DontPrintSemicolon%
\SetAlgoVlined%
$n \leftarrow 0$\;
$\widetilde V_n \leftarrow \spanset{\{0\}}$\;
\While{
  not converged
}%
{
  \tcc{solve reduced system}
  \Find{$\widetilde u_n \in \widetilde V_n \st$}{
    $a(\widetilde u_n, \varphi) = f(\varphi) \qquad \forall \varphi \in \widetilde V_n$}
  \tcc{form residual}
  $R_n(\cdot) \leftarrow f(\cdot) - a(\widetilde u_n, \cdot)$\;
  \tcc{find maximum local residual}
  $k \leftarrow \argmax\limits_{i=1, \dots, N_D} \norm{R_n}_{O_i'}$\;
  \tcc{solve local enrichment problem}
  \Find{$\localenrichment \in O_k \st$}{
    $a(\localenrichment, \varphi) = R_n(\varphi) \qquad \forall \varphi \in O_k$}
  \tcc{form enriched space}
  $\widetilde V_{n+1} \leftarrow \widetilde V_n \oplus \spanset{\{\localenrichment\}}$\;
  $n \leftarrow n+1$\;
}
\caption{Residual Based Online Enrichment}
\label{algo:online_enrichment}
\end{algorithm2e}
\DecMargin{1em}

Starting with the nullspace $\widetilde V_0$, we construct
a sequence of subspaces of $V$ which we denote by
$\widetilde V_n$.
The full problem is reduced by Galerkin projection on these
reduced spaces.
We denote the solutions of the reduced problems by $\widetilde u_n$.
Each reduced space $V_n$ is constructed by enriching the previous
reduced space with an additional basis function $\localenrichment$, 
which lies in one of the localized spaces $O_i$.

By using an overlapping domain decomposition,
all reduced local spaces and local spaces
are subspaces of $V$ (assuming an extension with
zero to the whole domain). 
A non overlapping domain decomposition would
require local spaces with non-zero boundary
conditions which are not subspaces of $V$.
This would require a completely different 
treatment.
\subsection*{Residual Based Enrichment}
The residual based enrichment algorithm 
(given as Algorithm \ref{algo:online_enrichment})
first selects the local enrichment space from which the 
enrichment function $\localenrichment$ is taken.
The local space $O_k$ which maximizes the dual norm
of the residual $\norm{R_n}_{O_i'} = \sup_{\varphi \in O_i \setminus \{0\}} \frac{R_n(\varphi)}{\norm{\varphi}_a}$ is chosen,
i.e.
\begin{equation}
k:= \argmax\limits_{i=1, \dots, N_D} \norm{R_n}_{O_i'}.
\end{equation}
The residual $R_n \in V'$ is defined as
$R_n(\cdot) := f(\cdot) - a(\widetilde u_n, \cdot)$.
Then a localized problem is formed by 
a Galerkin projection
of the original problem onto this local space
$O_k$,
and replacing the right hand side $f$ by the last residual $R_n$.
The solution of the localized problem is the enrichment function $\localenrichment$.
\subsection*{Globally Coupled Local Enrichment}
The globally coupled local enrichment algorithm
(given as Algorithm \ref{algo:globally_coupled_online_enrichment})
couples the global reduced space with the full local space.
First it iterates over all local spaces $O_i$
and solves the coupled problem: It solves
the original problem projected on the space
$\widetilde V_n \oplus O_i$, the solution of
this coupled problem is called $\localenrichedsolution{i}$.
Then the local space $O_k$ is selected which maximizes
the change in the solution $\norm{\widetilde u_n - \localenrichedsolution{k}}_a$, i.e.~
\begin{equation}
k:= \argmax\limits_{i=1, \dots, N_D}
\norm{\widetilde u_n - \localenrichedsolution{i}}_a
.
\end{equation}
The function $\localenrichedsolution{k}$ is used to enrich the space $\widetilde V_n$.
Note that this is an enrichment in $O_k$, even though $\localenrichedsolution{k}$ has global support.

\subsection*{Runtimes}
While Algorithm \ref{algo:online_enrichment}
has to solve only one local problem in each iteration,
Algorithm \ref{algo:globally_coupled_online_enrichment}
has to solve $N_D$ coupled problems in each iteration.
However, depending on the context, 
Algorithm \ref{algo:globally_coupled_online_enrichment}
might be preferable, because the coupled problems
are still small (approximately of dimension
$\mathrm{dim}(\widetilde V_n) + \mathrm{dim}(O_i)$)
and they are all independent and can thus be solved in parallel.
\IncMargin{1em}
\begin{algorithm2e}
\DontPrintSemicolon%
\SetAlgoVlined%
$n \leftarrow 0$\;
$\widetilde V_n \leftarrow \spanset{\{0\}}$\;
\While{
  not converged
}%
{
  \tcc{solve reduced system}
  \Find{$\widetilde u_n \in \widetilde V_n \ \st$}{
    $a(\widetilde u_n, \varphi) = f(\varphi) \qquad \forall \varphi \in \widetilde V_n$}
  \tcc{solve local enriched problems}
  \For{$i = 1, \dots, N_D$}{
    \Find{$\localenrichedsolution{i} \in \widetilde V_n \oplus O_i \ \st$}{
      $a(\localenrichedsolution{i}, \varphi) = f(\varphi) \qquad \forall \varphi \in \widetilde V_n \oplus O_i$}
  }
  \tcc{find maximum solution shift}
  $k \leftarrow \argmax\limits_{i=1, \dots, N_D} \norm{\widetilde u_n - \localenrichedsolution{i}}_a$\;
  \tcc{form enriched space}
  $\widetilde V_{n+1} \leftarrow \widetilde V_n \oplus \spanset{\{\localenrichedsolution{k}\}}$\;
  $n \leftarrow n+1$\;
}
\caption{Globally Coupled Online Enrichment}
\label{algo:globally_coupled_online_enrichment}
\end{algorithm2e}
\DecMargin{1em}

\section{A priori convergence}
First we prove exponential convergence for the residual based enrichment.
\begin{thm}[Exponential convergence]
\label{maintheorem}
For the reduced solutions $\widetilde u_{n+1}$ in Algorithm \ref{algo:online_enrichment} it holds that
\begin{equation}
\norm{u - \widetilde u_{n+1}}_a
\leq c \cdot 
\norm{u - \widetilde u_{n}}_a
\end{equation}
with
\begin{equation}
\label{defc}
c := \sqrt{1 - \frac{1}{N_D} \frac{1}{c_{pu}^2}} .
\end{equation}
The constant $c_{pu}$ is explained and defined
later in this section.
\end{thm}

\begin{pf}
As we use the energy norm, the solution is the best approximation
\begin{equation}
\norm{\widetilde u_{n+1} - u}_a \leq \norm{\varphi - u}_a \qquad \forall \varphi \in \widetilde V_{n+1} .
\end{equation}
This holds that for $\varphi = \widetilde u_n + \alpha \localenrichment$ for all $\alpha$ in $\mathbb{R}$.
Because of the symmetry of $a$ it holds that
\begin{equation}
\nonumber
\norm{\widetilde u_n + \alpha \localenrichment - u}_a^2 = \norm{\widetilde u_n - u}_a^2 - 2 \alpha R_n(\localenrichment) + \alpha^2 \norm{\localenrichment}_a^2 .
\end{equation}
This term is minimized by choosing $\alpha = R_n(\localenrichment) / \norm{\localenrichment}_a^2$. We use this $\alpha$ and 
realize that
$R_n(\localenrichment) = \norm{\localenrichment}_a^2 = \norm{R_n}_{O_k'}^2$,
because $\localenrichment$ is the Riesz representative of $R_n$ in $O_k$ in the energy norm.
It follows that
\begin{equation}
\label{eq:chungend}
\norm{\widetilde u_{n+1} - u}_a^2 \leq \norm{\widetilde u_n - u}_a^2 - \norm{R_n}_{O_k'}^2 .
\end{equation}
Till this point, the proof followed the structure given in \cite[Section 4]{Chung2015a}. 
We defined $k$ to select the largest local residual, so it holds that
\begin{equation}
\label{eq:r1}
\norm{R_n}_{O_k'}^2 \geq \frac{1}{N_D} \sum_{i=1}^{N_D} \norm{R_n}_{O_i'}^2 .
\end{equation}
Furthermore, from \cite[Proposition 5.1]{Buhr2017} we know
\begin{equation}
\label{eq:r2}
\norm{R_n}_{V'}^2 \leq c_{pu}^2 \sum_{i=1}^{N_D} \norm{R_n}_{O_i'}^2 
\end{equation}
with a constant $c_{pu}$ which is a stability constant for the
partition of unity (pu). Magnitude and scaling behavior
of $c_{pu}$ depend on the choice of the partition of unity
and the norm of the spaces. See \cite{Buhr2017} for more details.
As we use the energy norm, we have 
\begin{equation}
\label{eq:r3}
\norm{R_n}_{V'}^2 = \norm{\widetilde u_n - u}_a^2 .
\end{equation}
Combining \eqref{eq:r1}, \eqref{eq:r2}, and \eqref{eq:r3} we obtain
\begin{equation}
\label{eq:rall}
\norm{R_n}_{O_k'}^2 \geq \frac{1}{N_D} \frac{1}{c_{pu}^2} \norm{\widetilde u_n - u}_a^2 .
\end{equation}
Combining \eqref{eq:rall} with \eqref{eq:chungend} yields
\begin{equation}
\norm{\widetilde u_{n+1} - u}_a^2 \leq \left( 1 - \frac{1}{N_D} \frac{1}{c_{pu}^2} \right) \norm{\widetilde u_n - u}_a^2
\end{equation}
and thus the claim. \hfill $\square$
\end{pf}

\begin{cor}
\label{cor1}
For the reduced solutions $\widetilde u_n$ in Algorithm \ref{algo:online_enrichment} it holds that
\begin{equation}
\norm{\widetilde u_{n} - u}_a
\leq c^n \cdot 
\norm{u}_a
\end{equation}
with $c$ as defined in \eqref{defc}.
\end{cor}
\begin{pf}
This follows from Theorem \ref{maintheorem}, because 
$\widetilde u_0 = 0$.\hfill $\square$

\end{pf}
\subsection*{Concerning $c_{pu}$}
\label{sec:concerningcpu}
The constant $c_{pu}$ is defined to be 
\begin{equation}
\label{eq:cpudef}
c_{pu}^2 := \sup_{\varphi \in V \setminus \{0\}} 
\frac{
  \sum_{i=1}^{N_D} \norm{\varrho_i \varphi}_a^2
}{
  \norm{\varphi}_a^2
}.
\end{equation}
With this constant, it holds that
\begin{equation}
\norm{\zeta}_{V'}^2 \leq c_{pu}^2 \sum_{i=1}^{N_D} \norm{\zeta}_{O_i'}^2
\end{equation}
for any element $\zeta$ of the dual space $V'$, especially for the residual $R_n$
(see \cite[Proposition 5.1]{Buhr2017} for details).
An upper bound for $c_{pu}$ is devised in the following proposition.
\begin{prop}[Upper bound of $c_{pu}$]
With $J$ being the maximum number of functions $\varrho_i$ having support
in any given point $x$ in $\Omega$, $\frac{\kappa_{max}}{\kappa_{min}}$
being the contrast of the problem, $c_f$ the constant of the
Friedrich's inequality on $\Omega$ and $\norm{\cdot}_\infty$
the infinity norm, it holds that
\begin{equation}
c_{pu}^2 \leq
2 J \left(c_f
\frac{\kappa_{max}}{\kappa_{min}}
\max_i \norm{\nabla \varrho_i}_\infty^2 + \max_i \norm{\varrho_i}_\infty^2
\right).
\end{equation}
\end{prop}
\begin{pf}
Starting from \eqref{eq:cpudef}, we estimate 
$\sum_{i=1}^{N_D} \norm{\varrho_i \varphi}_a^2$. It holds that
\begin{equation}
\sum_{i=1}^{N_D} \norm{\varrho_i \varphi}_a^2
=
\sum_{i=1}^{N_D} \int_{\omega_i} \kappa |\nabla(\varrho_i \varphi)|^2 dx
\end{equation}
\begin{equation}
\label{eqlast}
\leq
\sum_{i=1}^{N_D} \int_{\omega_i} \kappa \left[2|(\nabla \varrho_i) \varphi|^2 + 2|\varrho_i (\nabla \varphi)|^2 \right]dx .
\end{equation}
For the second term in \eqref{eqlast} it holds that
\begin{equation}
\label{part1}
\sum_{i=1}^{N_D} \int_{\omega_i} 2 \kappa |\varrho_i (\nabla \varphi)|^2 dx
\leq 2 J \max_i \norm{\varrho_i}_\infty^2 \norm{\varphi}_a^2
\end{equation}
and for the first term we have
\begin{equation}
\sum_{i=1}^{N_D} \int_{\omega_i} 2 \kappa |(\nabla \varrho_i) \varphi|^2 dx
\leq
2 J \kappa_{max} \max_i \norm{\nabla \varrho_i}_\infty^2 \int_\Omega |\varphi|^2 dx
\nonumber
\end{equation}
\begin{equation}
\label{part2}
\leq
2 J c_f \frac{\kappa_{max}}{\kappa_{min}} \max_i \norm{\nabla \varrho_i}_\infty^2 \norm{\varphi}_a^2
.
\end{equation}
Combining these yields the claim.
\hfill $\square$
\end{pf}
\subsection*{Globally coupled enrichment}
The globally coupled enrichment given in Algorithm \ref{algo:globally_coupled_online_enrichment}
is the optimal enrichment: Among all enrichment functions
from all local spaces, it selects the one which minimizes the
resulting error in the energy norm.
\begin{thm}[Optimality of Algorithm \ref{algo:globally_coupled_online_enrichment}]
For the reduced solutions $\widetilde u_{n+1}$ in Algorithm \ref{algo:globally_coupled_online_enrichment} it holds that
\begin{equation}
\begin{multlined}
\norm{u - \widetilde u_{n+1}}_a = 
\min_{i = 1, \dots, N_D}
\min_{\psi_e \in O_i}
\hfill \\
\Big\{
\norm{u - \widetilde u_e}_a
\ \Big| \ 
\widetilde u_e \in \widetilde V_n \oplus \spanset \{ \psi_e \}
\ \mathrm{solves} \hfill \\ 
\hspace{20pt} a(\widetilde u_e, \varphi) = f(\varphi) \qquad \forall \varphi \in \widetilde V_n  \oplus \spanset \{ \psi_e \}
\Big\} .
\end{multlined}
\end{equation}
\end{thm}
\begin{pf}
First we realize that the solution $\widetilde u_{n+1}$ is identical to 
$\localenrichedsolution{k}$, because $\localenrichedsolution{k}$ solves
$a(\localenrichedsolution{k}, \varphi) = f(\varphi)$ in
$\widetilde V_{n+1}$, which is a subspace of $\widetilde V_n \oplus O_k$ and the solution is unique:
\begin{equation}
\widetilde u_{n+1} = \localenrichedsolution{k} .
\end{equation}
Second, since $u - \localenrichedsolution{i}$ is $a$-orthogonal to
$\localenrichedsolution{i} - \widetilde u_n$, it holds that
\begin{equation}
\norm{u - \widetilde u_n}_a^2 = \norm{u - \localenrichedsolution{i}}_a^2 + \norm{\localenrichedsolution{i} - \widetilde u_n}_a^2
\end{equation}
and thus 
\begin{equation}
k = \argmax\limits_{i=1, \dots, N_D} \norm{\widetilde u_n 
- \localenrichedsolution{i}}_a
\end{equation}
implies
\begin{equation}
k = \argmin\limits_{i=1, \dots, N_D} \norm{u - \localenrichedsolution{i}}_a.
\end{equation}
So $\localenrichedsolution{k}$ is closest to $u$ among all $\localenrichedsolution{i}$.

Third, $\localenrichedsolution{i}$ is the best approximation in $\widetilde V_n \oplus O_i$.
It is not possible to get closer to $u$ with any other enrichment in $O_i$.
 \hfill $\square$
\end{pf}

\begin{cor}
The results in Theorem \ref{maintheorem} and Corollary \ref{cor1} also hold for Algorithm \ref{algo:globally_coupled_online_enrichment}.
\end{cor}
\begin{pf}
The enrichment in Algorithm \ref{algo:globally_coupled_online_enrichment} is optimal,
so it is not worse than the enrichment of Algorithm \ref{algo:online_enrichment}.
Any bound for Algorithm \ref{algo:online_enrichment} holds also for 
Algorithm \ref{algo:globally_coupled_online_enrichment}.
 \hfill $\square$
\end{pf}
\section{Numerical Experiments}
\begin{figure}
\centering
{%
\setlength{\fboxsep}{0pt}%
\setlength{\fboxrule}{1pt}%
\fbox{\includegraphics[width=0.2\textwidth]{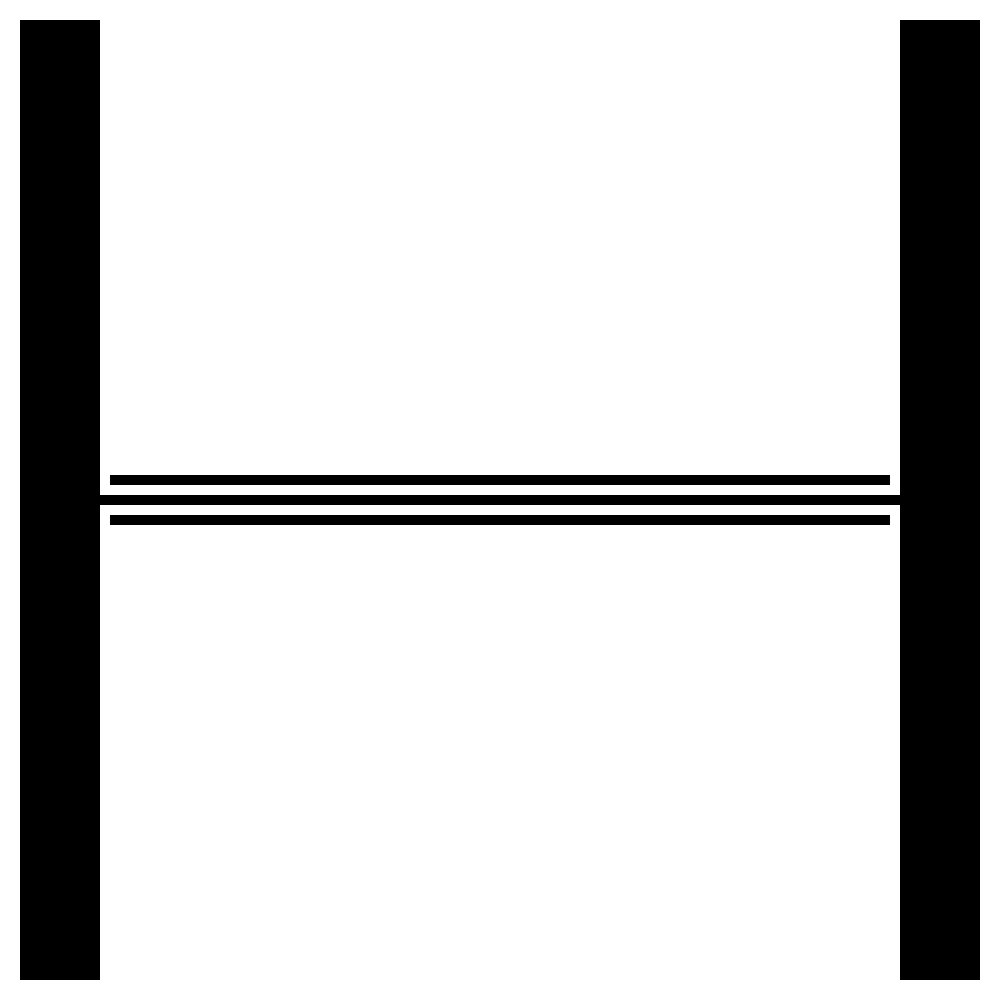}}
}%
\ \ \ 
\includegraphics[width=0.2\textwidth]{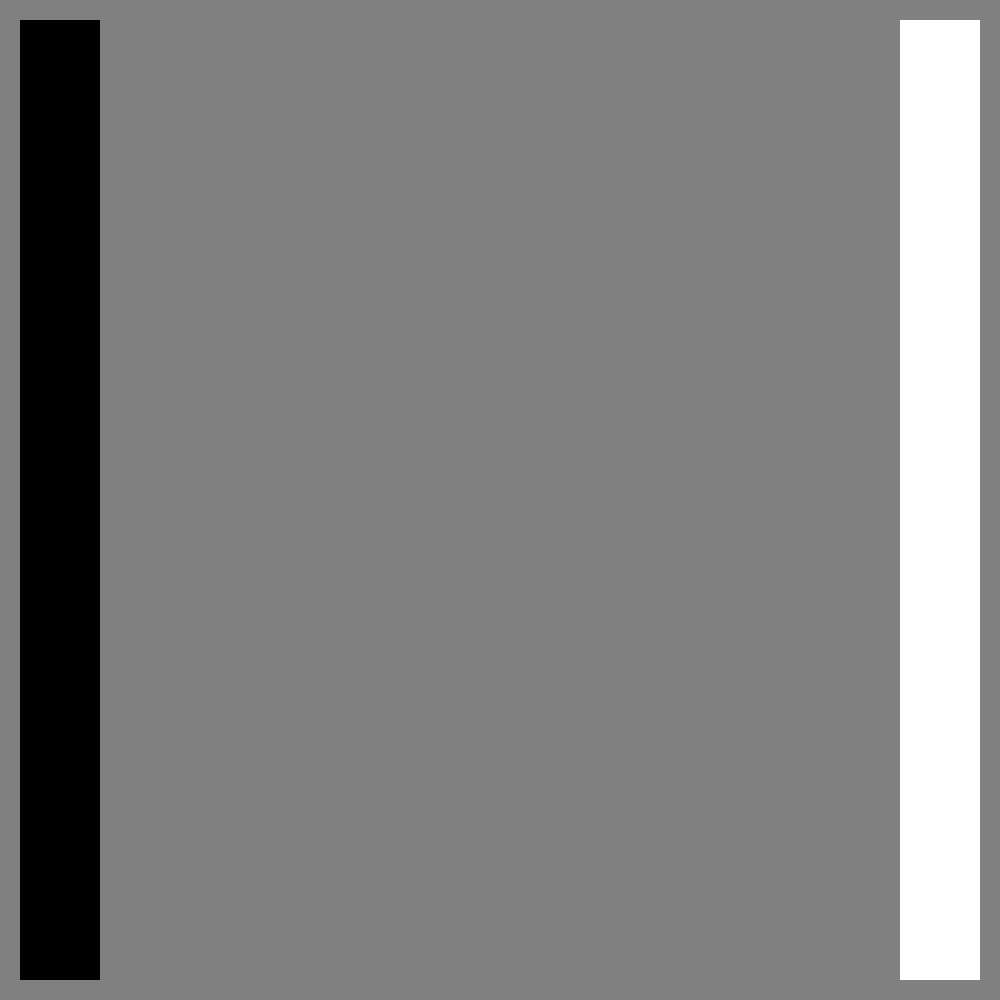}
\caption{Left: Coefficient field $\kappa$. White is 1, black is $10^5$. Right: right hand side $f$. Black is $-10^5$, gray is $0$, white is $10^5$.}
\label{fig:geometry}
\end{figure}

\begin{figure}
\centering
\includegraphics[width=0.34\textwidth]{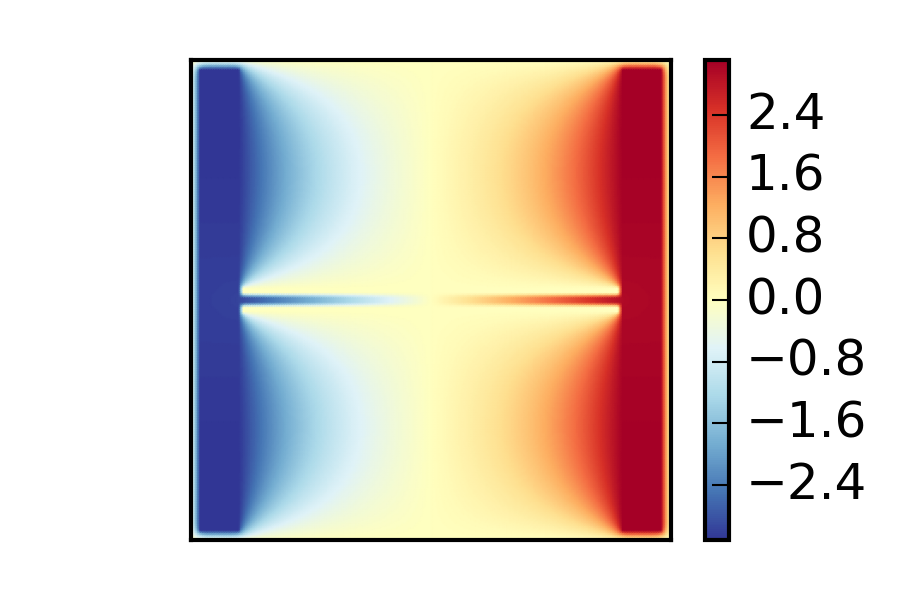}
\caption{Reference solution of problem.}
\label{fig:solution}
\end{figure}

\begin{figure}
\centering
\def\mysize{4}
\begin{tikzpicture}[scale=\mysize]
\node at (0.5,0.5) {\includegraphics[width=\mysize cm]{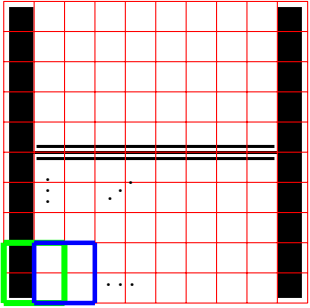}};
\draw [step=0.1, red] (0.,0.) grid (1,1) ;
\draw [step=0.2, green, line width=0.8mm] (0.,0.) grid (0.2,0.2) ;
\draw [step=0.2, blue, line width=0.6mm, shift={(0.1,0)}] (0.0,0.0) grid (0.2,0.2) ;
\node at (0.3, 0.1) [thick, anchor=north west] {$\dots$};
\node at (0.1, 0.3) [thick, anchor=south west] {$\vdots$};
\node at (0.3, 0.3) [thick, anchor=south west] {\reflectbox{$\ddots$}};
\end{tikzpicture}
\caption{Overlapping domain decomposition $\omega_i$ (green and blue) constructed from $2\times 2$ patches
of non-overlapping domain decomposition (red).}
\label{fig:dd}
\end{figure}

The experiments were carried out using the software package
pyMOR (\cite{Milk2016}). The source code to reproduce
the results presented here can be obtained at Zenodo (\cite{Buhr2017b}).

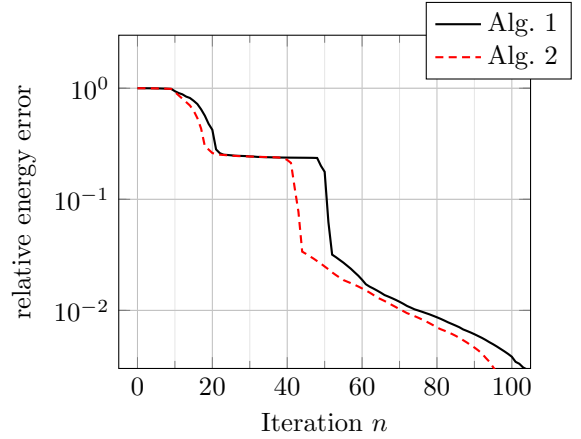
\begin{figure}
\centering
\begin{tikzpicture}
\begin{semilogyaxis}[
    width=7cm,
    height=6cm,
    xmin=-5,
    xmax=105,
    ymin=3e-3,
    ymax=3e0,
    xlabel=Iteration $n$,
    ylabel=relative energy error,
    legend pos=north east,
    grid=both,
    grid style={line width=.1pt, draw=gray!20},
    major grid style={line width=.2pt,draw=gray!50},
    ytick={1e-6, 1e-5, 1e-4, 1e-3, 1e-2, 1e-1, 1e0},
    minor ytick={1e-15, 1e-14, 1e-13, 1e-12, 1e-11, 1e-10, 1e-9, 1e-8, 1e-7, 1e-6, 1e-5, 1e-4, 1e-3, 1e-2, 1e-1, 1e0, 1e1, 1e2, 1e3},
    minor xtick={10,30,50,70,90},
    legend style={at={(1.1,1.1)},anchor=north east},
  ]
\addplot [solid, thick] table [x expr=\coordindex,y index=0] {errors.dat};
\addplot [densely dashed, thick, red] table [x expr=\coordindex,y index=0] {errors_g_c.dat};
\legend{Alg. \ref{algo:online_enrichment}, Alg. \ref{algo:globally_coupled_online_enrichment}}
\end{semilogyaxis}
\end{tikzpicture}
\caption{Decay of relative energy error during iteration (zoom).}
\label{fig:error decay zoom}
\end{figure}

\begin{figure}
\centering
\begin{tikzpicture}
\begin{semilogyaxis}[
    width=7cm,
    height=6cm,
    xmax=610,
    xlabel=Iteration $n$,
    ylabel=relative energy error,
    legend pos=north east,
    grid=both,
    grid style={line width=.1pt, draw=gray!20},
    major grid style={line width=.2pt,draw=gray!50},
    ytick={1e-12, 1e-9, 1e-6, 1e-3, 1},
    minor ytick={1e-15, 1e-14, 1e-13, 1e-12, 1e-11, 1e-10, 1e-9, 1e-8, 1e-7, 1e-6, 1e-5, 1e-4, 1e-3, 1e-2, 1e-1, 1e0, 1e1, 1e2, 1e3},
    xtick={0,100,200,300,400,500,600},
    minor xtick={50,150,250,350,450,550,650},
    legend style={at={(1.1,1.1)},anchor=north east},
  ]
\addplot [solid, thick] table [x expr=\coordindex,y index=0] {errors.dat};
\addplot [densely dashed, thick, red] table [x expr=\coordindex,y index=0] {errors_g_c.dat};
\draw [red, thick,rounded corners] (axis cs:-5,3e-3) rectangle (axis cs:105,3e0);
\node at (axis cs:120,1e-0) [red, anchor=north west] {Fig.~\ref{fig:error decay zoom}};
\fill [red, opacity=0.2] (420,1e-20) rectangle (1000, 1000);
\node at (430, 1e-9) [red, anchor=north west, rotate=90] {numerical noise};
\legend{Alg. \ref{algo:online_enrichment}, Alg. \ref{algo:globally_coupled_online_enrichment}}
\end{semilogyaxis}
\end{tikzpicture}
\caption{Decay of relative energy error during iteration.}
\label{fig:error decay}
\end{figure}
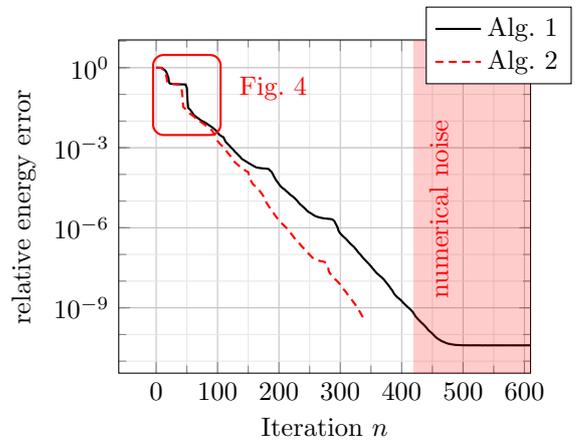

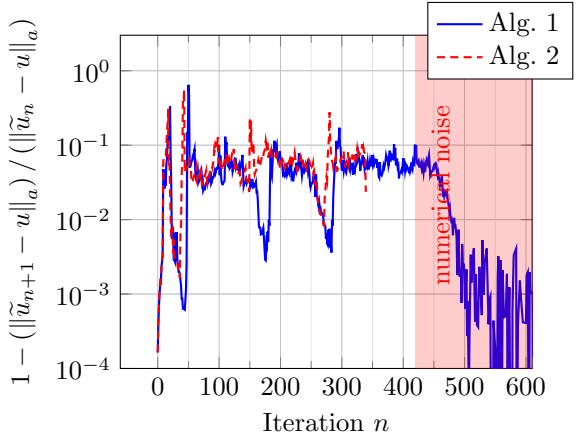
\begin{figure}
\centering
\begin{tikzpicture}
\begin{semilogyaxis}[
    width=7cm,
    height=6cm,
    xmax=610,
    ymin = 1e-4,
    ymax=3e0,
    xlabel=Iteration $n$,
    ylabel=$1 - \left(\norm{\widetilde u_{n+1} - u}_a\right)/\left(\norm{\widetilde u_n - u}_a\right)$,,
    legend pos=north east,
    grid=both,
    grid style={line width=.1pt, draw=gray!20},
    major grid style={line width=.2pt,draw=gray!50},
    ytick={1e-15, 1e-14, 1e-13, 1e-12, 1e-11, 1e-10, 1e-9, 1e-8, 1e-7, 1e-6, 1e-5, 1e-4, 1e-3, 1e-2, 1e-1, 1e0, 1e1, 1e2, 1e3},
    minor ytick={1e-15, 1e-14, 1e-13, 1e-12, 1e-11, 1e-10, 1e-9, 1e-8, 1e-7, 1e-6, 1e-5, 1e-4, 1e-3, 1e-2, 1e-1, 1e0, 1e1, 1e2, 1e3},
    xtick={0,100,200,300,400,500,600},
    minor xtick={50,150,250,350,450,550,650},
    legend style={at={(1.1,1.1)},anchor=north east},
  ]
\addplot [solid, thick, blue] table [x expr=\coordindex,y index=0] {convergence.dat};
\addplot [densely dashed, thick, red] table [x expr=\coordindex,y index=0] {convergence_g_c.dat};
\draw [red, thick] (axis cs:-100,1.714e-10) -- (axis cs:610,1.714e-10); %
\node at (axis cs:-40,1e-8) [red, anchor=north west] {guaranteed rate $(1-c)$};
\fill [red, opacity=0.2] (420,1e-20) rectangle (1000, 1000);
\node at (430, 1e-3) [red, anchor=north west, rotate=90] {numerical noise};
\legend{Alg. \ref{algo:online_enrichment}, Alg. \ref{algo:globally_coupled_online_enrichment}}
\end{semilogyaxis}
\end{tikzpicture}
\caption{Convergence. 1 would be convergence within one iteration, 0 would be stagnation.}
\label{fig:convergence rate}
\end{figure}

While the method and the proofs work both in two and three dimensions,
experiments were conducted for the two dimensional case only.
We define $\Omega$ to be the unit square $(0,1)^2$.
We discretize the problem using P1 finite elements on 
a regular grid of $200\times 200$ squares, each divided
into four triangles, resulting in 80,401 degrees of freedom.
We use a coefficient field $\kappa$ with high contrast
($\kappa_{max}/\kappa_{min} = 10^5$) and high
conductivity channels to get interesting behavior
(see Fig.~\ref{fig:geometry} and \ref{fig:solution}).
As domain decomposition $\omega_i$ we use domains of size 
$0.2\times 0.2$ with overlap $0.1$, resulting in $81$
subdomains (Fig.~\ref{fig:dd}).
The resulting error decay is shown in Fig.~\ref{fig:error decay zoom}, Fig.~\ref{fig:error decay}, and \ref{fig:convergence rate}.
To compare with the theory,
we calculate an upper bound for $c_{pu}$ using
$J=4$,
$c_f = 1/(\sqrt{2}\pi)$,
$\kappa_{max} / \kappa_{min} = 10^5$,
$\max_i \norm{\nabla \varrho_i}_\infty^2 = 2 H^{-2} = 200$,
$\max_i \norm{\varrho_i}_\infty^2 = 1$
and obtain $c_{pu}^2 \leq 3.6013 \cdot 10 ^ 7$.
This results in an estimate of $1-c \geq 1.714 \cdot 10^{-10}$.
The rate of convergence observed in the experiment is several
orders of magnitude better than the rate guaranteed by the
a priori theory
and is close to the optimal convergence rate
(Fig.~\ref{fig:convergence rate}).

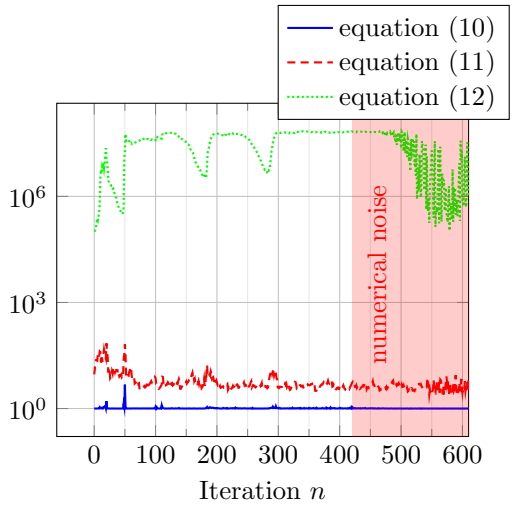
\begin{figure}
\centering
\begin{tikzpicture}
\begin{semilogyaxis}[
    width=7cm,
    height=6cm,
    xmax=610,
    xlabel=Iteration $n$,
    legend pos=north east,
    grid=both,
    grid style={line width=.1pt, draw=gray!20},
    major grid style={line width=.2pt,draw=gray!50},
    xtick={0,100,200,300,400,500,600},
    minor xtick={50,150,250,350,450,550,650},
    legend style={at={(1.1,0.95)},anchor=south east},
  ]
\addplot [solid, thick, blue] table [x expr=\coordindex,y index=0] {ineq.dat};
\addplot [densely dashed, thick, red,] table [x expr=\coordindex,y index=1] {ineq.dat};
\addplot [densely dotted, thick, green,] table [x expr=\coordindex,y index=2] {ineq.dat};
\fill [red, opacity=0.2] (420,1e-20) rectangle (1000, 1e9);
\node at (430, 1e1) [red, anchor=north west, rotate=90] {numerical noise};
\legend{
  equation \eqref{eq:chungend},
  equation \eqref{eq:r1},
  equation \eqref{eq:r2},
}
\end{semilogyaxis}
\end{tikzpicture}
\caption{Sharpness of inequalities in Algorithm \ref{algo:online_enrichment}:
For   equation \eqref{eq:chungend},
  $\left(\norm{\widetilde u_n - u}_a^2 - \norm{R_n}_{O_k'}^2\right)/\left(
  \norm{\widetilde u_{n+1} - u}_a^2\right)$ is plotted.
For   equation \eqref{eq:r1},
$\left(\norm{R_n}_{O_k'}^2\right) / \left(
\frac{1}{N_D} \sum_{i=1}^{N_D} \norm{R_n}_{O_k'}^2\right) $ is plotted.
For   equation \eqref{eq:r2},
$\left(c_{pu}^2 \sum_{i=1}^{N_D} \norm{R_n}_{O_i'}^2\right)/\left(
\norm{R_n}_{V'}^2\right)$ is plotted.
}
\label{fig:ineq}
\end{figure}

To investigate the reason for this,
we plot the quotient of the larger part and the smaller part
of the estimates \eqref{eq:chungend}, \eqref{eq:r1},
and \eqref{eq:r2} in Fig.~\ref{fig:ineq}.
It can be observed that the estimate \eqref{eq:chungend}
is rather sharp, except when the error drops after
a plateau.
In estimate \eqref{eq:r1}, around one order of magnitude is lost.
This could be improved by not enriching only one space but
using a marking strategy instead.
However, the main reason for the a priori theory to be so 
pessimistic seems to be in estimate \eqref{eq:r2}.

\section{Conclusion and Outlook}
We have shown that residual based
online enrichment %
converges exponentially.
The observed rate of convergence is far better than
the rate guaranteed by theory and close to the
rate of optimal convergence.
However, these results do not transfer immediately to
methods like ArbiLoMod or LRBMS,
because these methods do not enrich with a local solution,
but they apply a subspace projection to the local solution
before adding it to the reduced basis.
The convergence behavior with an additional subspace projection
is subject to future work.
Additionally, the results presented in this publication
are restricted to not parameterized problems. Also the extension to 
parameterized problems is an interesting question which
remains to be answered.
\bibliography{/home/andreasbuhr/bibliothek/refs}

\begin{thebibliography}{10}
\providecommand{\natexlab}[1]{#1}
\providecommand{\url}[1]{\texttt{#1}}
\providecommand{\urlprefix}{URL }
\expandafter\ifx\csname urlstyle\endcsname\relax
  \providecommand{\doi}[1]{doi:\discretionary{}{}{}#1}\else
  \providecommand{\doi}{doi:\discretionary{}{}{}\begingroup
  \urlstyle{rm}\Url}\fi

\bibitem[{Albrecht et~al.(2012)Albrecht, Haasdonk, Ohlberger, and
  Kaulmann}]{Albrecht2012}
Albrecht, F., Haasdonk, B., Ohlberger, M., and Kaulmann, S. (2012).
\newblock The localized reduced basis multiscale method.
\newblock \emph{Proceedings of Algoritmy 2012, Conference on Scientific
  Computing, Vysoke Tatry, Podbanske, September 9-14, 2012}, 393--403.

\bibitem[{Albrecht and Ohlberger(2013)}]{Albrecht2013}
Albrecht, F. and Ohlberger, M. (2013).
\newblock The localized reduced basis multi-scale method with online
  enrichment.
\newblock \emph{Oberwolfach Rep.}, 7, 406--409.
\newblock \doi{10.4171/OWR/2013/07}.

\bibitem[{Buhr(2017)}]{Buhr2017b}
Buhr, A. (2017).
\newblock {Source Code to "Exponential Convergence of Online Enrichment in
  Localized Reduced Basis Methods"}.
\newblock \doi{10.5281/zenodo.1002767}.

\bibitem[{Buhr et~al.(2017)Buhr, Engwer, Ohlberger, and Rave}]{Buhr2017}
Buhr, A., Engwer, C., Ohlberger, M., and Rave, S. (2017).
\newblock {ArbiLoMod, a Simulation Technique Designed for Arbitrary Local
  Modifications}.
\newblock \emph{SIAM J. Sci. Comput.}, 39(4), A1435--A1465.
\newblock \doi{10.1137/15M1054213}.

\bibitem[{Chung et~al.(2015)Chung, Efendiev, and Leung}]{Chung2015a}
Chung, E.T., Efendiev, Y., and Leung, W.T. (2015).
\newblock Residual-driven online generalized multiscale finite element methods.
\newblock \emph{J. Comput. Phys.}, 302, 176--190.
\newblock \doi{10.1016/j.jcp.2015.07.068}.

\bibitem[{Chung et~al.(2017{\natexlab{a}})Chung, Efendiev, and
  Leung}]{Chung2017a}
Chung, E.T., Efendiev, Y., and Leung, W.T. (2017{\natexlab{a}}).
\newblock Constraint energy minimizing generalized multiscale finite element
  method.
\newblock \emph{arXiv preprint arXiv:1704.03193}.

\bibitem[{Chung et~al.(2017{\natexlab{b}})Chung, Efendiev, and
  Leung}]{Chung2017c}
Chung, E.T., Efendiev, Y., and Leung, W.T. (2017{\natexlab{b}}).
\newblock Fast online generalized multiscale finite element method using
  constraint energy minimization.
\newblock \emph{arXiv preprint arXiv:1706.07093}.

\bibitem[{Efendiev et~al.(2013)Efendiev, Galvis, and Hou}]{Efendiev2013}
Efendiev, Y., Galvis, J., and Hou, T.Y. (2013).
\newblock Generalized multiscale finite element methods ({GMsFEM}).
\newblock \emph{J. Comput. Phys.}, 251, 116--135.
\newblock \doi{10.1016/j.jcp.2013.04.045}.

\bibitem[{Milk et~al.(2016)Milk, Rave, and Schindler}]{Milk2016}
Milk, R., Rave, S., and Schindler, F. (2016).
\newblock py{MOR} -- {G}eneric {A}lgorithms and {I}nterfaces for {M}odel
  {O}rder {R}eduction.
\newblock \emph{SIAM J. Sci. Comput.}, 38(5), S194--S216.
\newblock \doi{10.1137/15M1026614}.

\bibitem[{Ohlberger and Schindler(2015)}]{Ohlberger2015}
Ohlberger, M. and Schindler, F. (2015).
\newblock Error control for the localized reduced basis multi-scale method with
  adaptive on-line enrichment.
\newblock \emph{SIAM J. Sci. Comput.}, 37(6), A2865--A2895.
\newblock \doi{10.1137/151003660}.

\end{thebibliography}
\end{document}